\documentclass[12pt,twoside]{article}
\usepackage{times}
\usepackage{amsmath,amssymb}
\usepackage{amsthm}
\usepackage{mathrsfs}
\usepackage{color}
\usepackage{exscale}
\usepackage{relsize}
\usepackage{color}
\usepackage{fancyhdr}
\allowdisplaybreaks

\newcommand{\pf}{\noindent {\bf Proof. \hspace{2mm}}}

\newcommand{\be}{\begin{equation}}
\newcommand{\ee}{\end{equation}}
\newcommand{\bea}{\begin{eqnarray}}
\newcommand{\eea}{\end{eqnarray}}

\def\th{\theta}

\def\p{\partial}

 \allowdisplaybreaks

\begin{document}
 \footskip=0pt
 \footnotesep=2pt
\let\oldsection\section
\renewcommand\section{\setcounter{equation}{0}\oldsection}
\renewcommand\thesection{\arabic{section}}
\renewcommand\theequation{\thesection.\arabic{equation}}
\newtheorem{theorem}{\noindent Theorem}[section]
\newtheorem{lemma}{\noindent Lemma}[section]
\newtheorem*{claim}{\noindent Claim}
\newtheorem{proposition}{\noindent Proposition}[section]
\newtheorem*{inequality}{\noindent Weighted Inequality}
\newtheorem{definition}{\noindent Definition}[section]
\newtheorem{remark}{\noindent Remark}[section]
\newtheorem{corollary}{\noindent Corollary}[section]
\newtheorem{example}{\noindent Example}[section]

\title{On the Vanishing of some D-solutions to the Stationary Magnetohydrodynamics System}

\author{Zijin Li $^{a,b,}$\footnote{E-mail:zijinli@smail.nju.edu.cn} ,\quad Xinghong Pan$^{c,}$\footnote{E-mail:xinghong{\_}87@nuaa.edu.cn}\vspace{0.5cm}\\
 \footnotesize $^a$College of Mathematics and Statistics, Nanjing University of Information Science $\&$ Technology,\\
  \footnotesize Nanjing 210044, China.\\
  \footnotesize $^b$Department of Mathematics and IMS, Nanjing University, Nanjing 210093, China.\\
 \footnotesize $^c$Department of Mathematics, Nanjing University of Aeronautics and Astronautics, Nanjing 211106, China.\\
\vspace{0.5cm}
}

\date{}

\maketitle

\centerline {\bf Abstract}
In this paper, we study the stationary magnetohydrodynamics system in $\mathbb{R}^2\times\mathbb{T}$. We prove trivialness of D-solutions (the velocity field $u$ and the magnetic field $h$) when they are swirl-free. Meanwhile, this Liouville type theorem also holds provided $u$ is swirl-free and $h$ is axially symmetric, or both $u$ and $h$ are axially symmetric. Our method is also valid for certain related boundary value problems in the slab $\mathbb{R}^2\times[-\pi,\,\pi]$.
\vskip 0.3 true cm

\vskip 0.3 true cm

{\bf Keywords:} incompressible; magnetohydrodynamics system; swirl-free; axially symmetric.
\vskip 0.3 true cm

{\bf Mathematical Subject Classification 2010:} 35Q30, 76N10

\section{Introduction}
In this paper, we consider the stationary magnetohydrodynamics system
\begin{equation}\label{MHD}
\left\{
\begin{aligned}
&u\cdot\nabla u+\nabla p-h\cdot\nabla h-\Delta u=0,\\
&u\cdot\nabla h-h\cdot\nabla u-\Delta h=0,\\
&\nabla\cdot u=\nabla\cdot h=0,\\
\end{aligned}
\right.
\end{equation}
in $\mathbb{R}^2\times\mathbb{T}$ or in the slab $\mathbb{R}^2\times[-\pi,\,\pi]$, where $u(x),\,h(x)\in\mathbb{R}^3$, $p(x)\in\mathbb{R}$ represent the velocity vector, the magnetic field and the scalar pressure respectively. The MHD equations, which describe the state of the fluid flows of plasma, are fundamental partial differential equations in nature. For the background of the MHD system, we refer readers to \cite{Laudau1984} for more details. We note that if $h\equiv0$, the MHD system is reduced to the Navier-Stokes system.

In the following, we will carry out our proof in the cylindrical coordinates $(r,\,\th,\,z)$. That is, for $x=(x_1,\,x_2,\,x_3)\in\mathbb{R}^3$
\be
\left\{
\begin{aligned}
&r=\sqrt{x_1^2+x_2^2},\\
&\th=\arctan\frac{x_2}{x_1}, \\
&z=x_3.\\
\end{aligned}
\right.
\ee
And the solution of the incompressible stationary magnetohydrodynamics system is given as
\[
u=u^r(r,\th,z)e_r+u^{\th}(r,\th,z)e_{\th}+u^z(r,\th,z)e_z,
\]
\[
h=h^r(r,\th,z)e_r+h^{\th}(r,\th,z)e_{\th}+h^z(r,\th,z)e_z,
\]
where the basis vectors $e_r,e_\th,e_z$ are
\[
e_r=(\frac{x_1}{r},\frac{x_2}{r},0),\quad e_\th=(-\frac{x_2}{r},\frac{x_1}{r},0),\quad e_z=(0,0,1).
\]
The components $u^r,u^\th,u^z,h^r,h^\th,h^z$ satisfy
\begin{equation}\label{1.2}
\scriptsize
\left\{
\begin{aligned}
&(u^r\p_r+\frac{1}{r}u^\th\p_\th+u^z\p_z)u^r -\frac{(u^\th)^2}{r}+\frac{2}{r^2}\p_\th u^\th+\p_r p=(h^r\p_r+\frac{1}{r}h^\th\p_\th+h^z\p_z)h^r-\frac{(h^\th)^2}{r}+(\Delta-\frac{1}{r^2})u^r, \\
&(u^r\p_r+\frac{1}{r}u^\th\p_\th+u^z\p_z) u^\th+\frac{u^\th u^r}{r}-\frac{2}{r^2}\p_\th u^r+\frac{1}{r}\p_\th p=(h^r\p_r+\frac{1}{r}h^\th\p_\th+h^z\p_z)h^\th+\frac{h^rh^\th}{r}+(\Delta-\frac{1}{r^2})u^\th , \\
&(u^r\p_r+\frac{1}{r}u^\th\p_\th+u^z\p_z)u^z+\p_z p=(h^r\p_r+\frac{1}{r}h^\th\p_\th+h^z\p_z)h^z+\Delta u^z ,                                    \\
&(u^r\p_r+\frac{1}{r}u^\th\p_\th+u^z\p_z)h^r-(h^r\p_r+\frac{1}{r}h^\th\p_\th+h^z\p_z)u^r+\frac{2}{r^2}\p_\th h^\th=(\Delta-\frac{1}{r^2})h^r,\\
&(u^r\p_r+\frac{1}{r}u^\th\p_\th+u^z\p_z)h^\th-(h^r\p_r+\frac{1}{r}h^\th\p_\th+h^z\p_z)u^\th+\frac{u^\th h^r}{r}-\frac{h^\th u^r}{r}-\frac{2}{r^2}\p_\th h^r=(\Delta-\frac{1}{r^2})h^\th,\\
&(u^r\p_r+\frac{1}{r}u^\th\p_\th+u^z\p_z)h^z-(h^r\p_r+\frac{1}{r}h^\th\p_\th+h^z\p_z)u^z=\Delta h^z,\\
&\nabla\cdot u=\p_ru^r+\frac{u^r}{r}+\frac{1}{r}\p_\th u^\th+\p_zu^z=0,\quad \nabla\cdot h=\p_rh^r+\frac{h^r}{r}+\frac{1}{r}\p_\th h^\th+\p_zh^z=0.
\end{aligned}
\right.
\end{equation}
Here
\be
\Delta=\frac{\p^2}{\p r^2}+\frac{1}{r}\frac{\p}{\p r}+\frac{1}{r^2}\frac{\p^2}{\p\th^2}+\frac{\p^2}{\p z^2}
\ee
is the usual Laplacian operator.

The main aim of our paper is to study the Liouville type theorem of D-solutions of the stationary MHD system \eqref{1.2}. The study is partly motivated by the related Liouville problem of the stationary Navier-Stokes equations, which has attracted much attention in recent years and is still far from being fully understood. See for example \cite{CPZ2018, CPZ2018-1, Chae2014, Chae2016, Galdi2011, Koch2009, Korobkov2015, Kozono2017, Seregin2016, Seregin2018} and the reference therein. First, in full 3D case, the Liouville-type theorem holds provided the vanishing of $u^\th$ and $h^\th$. That is:
\begin{theorem}\label{thm1.1}
Let $(u,h)$ be a smooth solution to the problem
\begin{equation}\label{1.3}
\left\{
\begin{aligned}
&u\cdot\nabla u+\nabla p-h\cdot\nabla h-\Delta u=0,\quad\text{ in }\quad\mathbb{R}^2\times\mathbb{T}, \\
&u\cdot\nabla h-h\cdot\nabla u-\Delta h=0,\quad\text{ in }\quad\mathbb{R}^2\times\mathbb{T},\\
&\nabla\cdot u=\nabla\cdot h=0,\quad\text{ in }\quad\mathbb{R}^2\times\mathbb{T},\\
&u(x',z)=u(x',z+2\pi);\quad h(x',z)=h(x',z+2\pi),\\
&\lim_{|x|\to\infty}|u(x)|=0;\quad\lim_{|x|\to\infty}|h(x)|=0,\\
\end{aligned}
\right.
\end{equation}
with finite Dirichlet integral
\be\label{D-COND}
\int_\mathbb{T}\int_{\mathbb{R}^2}|\nabla u(x)|^2dx+\int_\mathbb{T}\int_{\mathbb{R}^2}|\nabla h(x)|^2dx<\infty.
\ee
Then $(u,h)\equiv0$ provided $u^\th=h^\th\equiv0$.
\end{theorem}

\qed
\begin{remark}
We emphasize here that our assumption of the smoothness of the solution $(u,\,h)$ is reasonable since one can derive the smoothness of any weak solution to \eqref{1.3} satisfying the D-condition \eqref{D-COND} by following the method developed in \cite{Galdi2011}.
\end{remark}

\qed

In the cylinder coordinate, we say a 3 dimensional vector field
\be
v(x)=v^r(r,\th,z)\cdot e_r+v^z(r,\th,z)\cdot e_z+v(r,\th,z)^\th e_\th
\ee
is axially symmetric if and only if
\be
\p_\th v^r=\p_\th v^z=\p_\th v^\th\equiv0.
\ee
Moreover, for axially symmetric magnetic field or axially symmetric velocity and magnetic fields, we derive two further results:
\begin{corollary}\label{Cor1}
Let $(u,h)$ be a smooth solution to the problem \eqref{1.3} with finite Dirichlet integral \eqref{D-COND}. Then $(u,h)\equiv0$ provided one of the following two conditions is satisfied:

(i) $u^\th\equiv0$ and $h$ is axially symmetric;

(ii) Both $u$ and $h$ are axially symmetric.
\end{corollary}

\qed

\begin{remark}
Consider the special case that $h\equiv0$, part (ii) of the above corollary is reduced to Theorem 1.1 in \cite{CPZ2018-1}. For part (i) ( also Theorem \ref{thm1.1} ) with $h\equiv0$, this Liouville-type theorem do not need to add the axially symmetric condition of $u$. More precisely, this is a result of the swirl-free full 3-D case.
\end{remark}

\qed

Instead of $u$ and $h$ are z-periodic, our method is valid for D-solutions of certain boundary value problems of magnetohydrodynamics system \eqref{MHD} in the slab $\mathbb{R}^2\times[-\pi,\,\pi]$. Here is the corollary:

\begin{corollary}\label{Cor2}
Let $(u,h)$ be a smooth solution to the magnetohydrodynamics system
\begin{equation}
\left\{
\begin{aligned}
&u\cdot\nabla u+\nabla p-h\cdot\nabla h-\Delta u=0,\quad\text{ in }\quad\mathbb{R}^2\times[-\pi,\pi], \\
&u\cdot\nabla h-h\cdot\nabla u-\Delta h=0,\quad\text{ in }\quad\mathbb{R}^2\times[-\pi,\pi],\\
&\nabla\cdot u=\nabla\cdot h=0,\quad\text{ in }\quad\mathbb{R}^2\times[-\pi,\pi],\\
&\lim_{|x|\to\infty}|u(x)|=0;\quad\lim_{|x|\to\infty}|h(x)|=0,\\
\end{aligned}
\right.
\end{equation}
with finite Dirichlet integral
\be
\int_{-\pi}^\pi\int_{\mathbb{R}^2}|\nabla u(x)|^2dx+\int_{-\pi}^\pi\int_{\mathbb{R}^2}|\nabla h(x)|^2dx<\infty
\ee
in the slab $\mathbb{R}^2\times[-\pi,\,\pi]$ equipped with the boundary conditions
\be
(u^z,\,\p_zu^r,\,\p_zu^\th)\Big|_{z\in\{-\pi,\,\pi\}}=0,\quad h\Big|_{z\in\{-\pi,\,\pi\}}=0, \\
\ee
or
\be
(u^z,\,\p_zu^r,\,\p_zu^\th)\Big|_{z\in\{-\pi,\,\pi\}}=0,\quad (h^z,\,\p_zh^r,\,\p_zh^\th)\Big|_{z\in\{-\pi,\,\pi\}}=0. \\
\ee
Then $(u,h)\equiv0$ provided one of the following three conditions is satisfied:

(i) $u^\th=h^\th\equiv0$;

(ii) $u^\th\equiv0$ and $h$ is axially symmetric;

(iii) Both $u$ and $h$ are axially symmetric.
\end{corollary}

\qed

We refer readers to the Appendix of our paper for some explanation to the reasonableness of the boundary conditions in the Corollary. The proof of Corollary \ref{Cor2} is omitted in this paper.

\begin{remark}
Unlike the magnetic field $h$, we do not know any vanishing result in a slab with $u$ satisfying homogeneous Dirichlet boundary condition.
\end{remark}

\qed

Our proof of the theorem and corollaries are based on the oscillation estimate of the pressure in \cite{CPZ2018-1}. Because of the partly "$\th-$dependent" of the pressure $p$ and some magnetic related terms, we need a careful treatment for getting the boundedness of $u$ and $h$ up to their second order derivatives and oscillation estimate of $p$ in a dyadic annulus. At last, we prove the Liouville type theorems by providing the vanishing of the $L^2$ norms of $\nabla u$ and $\nabla h$.


This paper is organized as follows. In Section 2, we give the proof of Theorem \ref{thm1.1}. Section 3 is devoted to proving the part (i) of Corollary \ref{Cor1}, while Section 4 is for the the part (ii). Some details of the boundary conditions in Corollary \ref{Cor2} could be found in the Appendix.

Throughout the paper, we use $C$ to denote a generic constant which may be different from line to line. We also apply $A\lesssim B$ to denote $A\leq CB$. We denote by $B(x_0,r):=\{x\in\mathbb{R}^d:|x-x_0|< r\}$. We simply denote by $B_r:=B(0,r)$ and $B:=B_1$.  For a domain $\Omega$ and $1\leq p\leq\infty$, $L^p(\Omega)$ denotes the usual Lebesgue space with norm $\|\cdot\|_{L^p(\Omega)}$. For $x=(x_1, x_2, x_3)\in\mathbb{R}^3$, we write $x=(x',x_3)$ or $(x',z)$ for simplicity. The symbol $\p_i$ stands for $\frac{\p}{\p x_i}$, for $i=1,2,3$, while $\p_r$, $\p_\th$ and $\p_z$ stands for $\frac{\p}{\p r}$, $\frac{\p}{\p \th}$ and $\frac{\p}{\p z}$ respectively.

\section{Proof of Theorem 1.1}
First we see, under the condition $u^\th=h^\th\equiv0$, \eqref{1.2} turns to
\begin{equation}\label{1.212}
\left\{
\begin{aligned}
&(u^r\p_r+u^z\p_z)u^r+\p_r p=(h^r\p_r+h^z\p_z)h^r+(\Delta-\frac{1}{r^2})u^r, \\
&-\frac{2}{r^2}\p_\th u^r+\frac{1}{r}\p_\th p=0, \\
&(u^r\p_r+u^z\p_z)u^z+\p_z p=(h^r\p_r+h^z\p_z)h^z+\Delta u^z ,                                    \\
&(u^r\p_r+u^z\p_z)h^r-(h^r\p_r+h^z\p_z)u^r=(\Delta-\frac{1}{r^2})h^r,\\
&(u^r\p_r+u^z\p_z)h^z-(h^r\p_r+h^z\p_z)u^z=\Delta h^z,\\
&\nabla\cdot u=\p_ru^r+\frac{u^r}{r}+\p_zu^z=0,\quad \nabla\cdot h=\p_rh^r+\frac{h^r}{r}+\p_zh^z=0.
\end{aligned}
\right.
\end{equation}
This section is divided into three parts. The first one is to derive the boundedness of $u$ and $h$ up to their second order derivatives. We have applied a result for local solutions in \cite{Zhang2006}. Second, by integrating the equation of $\p_r p$, we actually prove the boundedness of the oscillation of $p$ in a bounded dyadic annulus. Finally, by testing the MHD system with standard test functions, we prove the trivialness of $u$ and $h$.
\subsection{Boundedness of the solution up to second order derivatives}\label{SEC2.1}
\begin{lemma}\label{Lemma2.1}
Under the assumptions of Theorem \ref{thm1.1}, we have
\be
|\nabla^k u|+|\nabla^k h|\leq C_{k}<\infty,\quad 0\leq k\leq2.
\ee
Here $\nabla^kf$ denotes all the derivatives of $f$ with order $k$.
\end{lemma}
\pf Since $u$ and $h$ are assumed to be smooth functions and converge to 0 as $r\to\infty$, we have that both $u$ and $h$ are bounded. Now we derive the boundedness of their derivatives. By denoting
\be
w_1:=u+h;\quad w_2:=u-h,
\ee
\eqref{1.3} leads to
\be\label{EQW}
\left\{
\begin{aligned}
&w_2\cdot\nabla w_1-\Delta w_1+\nabla p=0,\\
&w_1\cdot\nabla w_2-\Delta w_2+\nabla p=0,\\
&\nabla\cdot w_1=\nabla\cdot w_2=0.
\end{aligned}
\right.
\ee
A direct application of Theorem 1.7 in \cite{Zhang2006} shows that, $\exists \,r_0\leq 1$, the gradient of $w_1$ and $w_2$ satisfy
\be
|\nabla w_i(x)|\leq\frac{C}{r_0^3}\int_{B(x,r_0)}|\nabla w_i(y)|dy+\frac{C}{r_0^4}\int_{B(x,r_0)}|w_i-(w_i)_{B(x,r_0)}|dy,\quad i=1,2.
\ee
Here, for $i=1,2$, $(w_i)_{B(x,r_0)}=\frac{1}{|B(x,r_0)|}\int_{B(x,r_0)}w_i(y)dy$. Applying Cauchy-Schwartz inequality and Poincar\'e inequality, one find
\be
\begin{split}
|\nabla w_i(x)|\lesssim&\frac{1}{r_0^3}\left(\int_{B(x,r_0)}|\nabla w_i(y)|^2dy\right)^{1/2}\left(\int_{B(x,r_0)}dy\right)^{1/2}\\
&+\frac{1}{r_0^4}\left(\int_{B(x,r_0)}|w_i-(w_i)_{B(x,r_0)}|^2dy\right)^{1/2}\left(\int_{B(x,r_0)}dy\right)^{1/2}\\
\lesssim&r_0^{-3/2}\|\nabla w_i\|_{L^2(B(x,1))}\\
\lesssim&r_0^{-3/2}(\|\nabla u\|_{L^2(B(x,1))}+\|\nabla h\|_{L^2(B(x,1))})\lesssim 1,\quad\text{ for }\quad i=1,2.
\end{split}
\ee
Taking the $curl$ of the first two equations of \eqref{EQW}, we then eliminate the terms of pressure. With the boundedness of $w_1$, $w_2$, $\nabla w_1$ and $\nabla w_2$, routine elliptic estimates prove the boundedness of $\nabla^2 w_1$ and $\nabla^2 w_2$. This leads to the boundedness of $u$ and $h$ up to their second order derivatives.

\qed

\subsection{Boundedness of the oscillation of $p$ in dyadic annulus}\label{sec2.3}
\begin{lemma}\label{Lemma2.2}
Under the assumptions of Theorem \ref{thm1.1}, for fixed $R>0$, it follows that
\be\
\sup_{r\in[R,\,2R],\,\th\in[0,\,2\pi],\,z\in[-\pi,\,\pi]}\left|p(r,\th,z)-p(R,0,0)\right|\lesssim1.
\ee
\end{lemma}

\pf First we note that the second equation of \eqref{1.212} turns to
\be
\p_\th p=\frac{2}{r}\p_\th u^r.
\ee
Owing to the boundedness of $\nabla u$, we have
\be
\left|\frac{1}{r}\p_\th u^r\right|=|-\sin\th\cdot\p_{x_1}u^r+\cos\th\cdot\p_{x_2}u^r|\leq|\nabla u|.
\ee
This leads to
\be
|\p_\th p|\lesssim 1.
\ee
Meanwhile, due to the third equation of \eqref{1.212} and Lemma \ref{Lemma2.1}, it follows
\be
\quad |\p_zp|\lesssim 1.
\ee
In the following part we will show that for any fixed $R>1$, the estimate
\be\label{2.2}
\left|\int_{-\pi}^\pi\int_0^{2\pi} p(r,\th,z)-p(R,\th,z)d\th dz\right|\lesssim 1
\ee
holds for all $r\in[R,\,2R]$. Now we integrate the first equation of \eqref{1.21} to get
\be
\begin{split}
\p_r\int_{-\pi}^\pi\int_0^{2\pi}pd\th dz=&\int_{-\pi}^\pi\int_0^{2\pi}\Big(-(u^r\p_r+u^z\p_z)u^r\Big.\\
&\hskip .5cm\Big.+(\p^2_r+\frac{1}{r}\p_r+\frac{1}{r^2}\p^2_\th+\p^2_z-\frac{1}{r^2})u^r+(h^r\p_r+h^z\p_z)h^r\Big) d\th dz\\
=&-\frac{1}{2}\int_{-\pi}^\pi\int_0^{2\pi}\p_r(u^r)^2d\th dz-\int_{-\pi}^\pi\int_0^{2\pi}u^z\p_zu^rd\th dz\\
&+\p_r^2\int_{-\pi}^\pi\int_0^{2\pi}u^rd\th dz+\frac{1}{r}\p_r\int_{-\pi}^\pi\int_0^{2\pi}u^rd\th dz\\
&-\frac{1}{r^2}\int_{-\pi}^\pi\int_0^{2\pi}u^rd\th dz+\frac{1}{2}\int_{-\pi}^\pi\int_{\mathbb{R}^2}\p_r(h^r)^2d\th dz\\
&+\int_{-\pi}^\pi\int_0^{2\pi}h^z\p_zh^rd\th dz.
\end{split}
\ee
$\forall$ $r_0\in[R,\, 2R]$, we integrate on $r$ from $R$ to $r_0$. It follows that
\be\label{2.4}
\begin{split}
&\left|\int_{-\pi}^\pi\int_0^{2\pi}p(r_0,\theta,z)-p(R,\th,z)d\th dz\right|\\
\lesssim&\left|\int_R^{r_0}\int_{-\pi}^\pi\int_0^{2\pi}\p_r(u^r)^2d\th dzdr\right|+\left|\int_R^{r_0}\int_{-\pi}^\pi\int_0^{2\pi}u^z\p_zu^rd\th dzdr\right|\\
&+\left|\int_R^{r_0}\p_r^2\int_{-\pi}^\pi\int_0^{2\pi}u^rd\th dzdr\right|+\left|\int_R^{r_0}\frac{1}{r}\p_r\int_{-\pi}^\pi\int_0^{2\pi}u^rd\th dzdr\right|\\
&+\left|\int_R^{r_0}\frac{1}{r^2}\int_{-\pi}^\pi\int_0^{2\pi}u^rd\th dzdr\right|+\left|\int_R^{r_0}\int_{-\pi}^\pi\int_{\mathbb{R}^2}\p_r(h^r)^2d\th dzdr\right|\\
&+\left|\int_R^{r_0}\int_{-\pi}^\pi\int_0^{2\pi}h^z\p_zh^rd\th dzdr\right|\\
:=&I_1+I_2+I_3+I_4+I_5+I_6+I_7.
\end{split}
\ee
In the following, we show that $I_1$ to $I_5$ in \eqref{2.4} are all bounded. First, due to the boundedness of $u$, we see
\be
I_1=\left|\int_{-\pi}^\pi\int_0^{2\pi}\left((u^r)^2(r_0,\th,z)-(u^r)^2(R,\th,z)\right)d\th dz\right|\lesssim 1.
\ee
Now we consider term $I_2$. Using integrating by parts and divergence free condition, we have
\be
\begin{split}
I_2=&\left|\int_R^{r_0}\int_{-\pi}^\pi\int_0^{2\pi}u^r\p_zu^zd\th dzdr\right|\\
=&\left|\int_R^{r_0}\int_{-\pi}^\pi\int_0^{2\pi}u^r\left(\p_ru^r+\frac{1}{r}u^r\right)d\th dzdr\right|\\
\lesssim&\,\,I_1+\left|\int_R^{r_0}\int_{-\pi}^\pi\int_0^{2\pi}\frac{(u^r)^2}{r}d\th dzdr\right|\\
\lesssim& 1+\left|\int_R^{2R}\frac{1}{r}dr\right|\lesssim 1.
\end{split}
\ee
Here the second inequality holds because the boundedness of $u$. For $I_3$, it follows
\be
\begin{split}
I_3=\left|\int_{-\pi}^\pi\int_0^{2\pi}\left(\p_ru^r(r_0,\th,z)-\p_ru^r(R,\th,z)\right)d\th dz\right|\lesssim1.
\end{split}
\ee
Meanwhile, $I_4$ satisfies the following estimate by using integration by parts
\be
\begin{split}
I_4&\lesssim\left|\int_{-\pi}^\pi\int_0^{2\pi}\frac{u^r}{r}\Big|^{r_0}_Rd\th dz\right|+\left|\int_R^{r_0}\int_{-\pi}^\pi\int_0^{2\pi}\frac{u^r}{r^2}d\th dzdr\right|\\
&\lesssim 1+\left|\int_R^{\infty}\frac{1}{r^2}dr\right|\lesssim 1.
\end{split}
\ee
Here the last two inequalities hold since $u^r$ is bounded. And the related estimate for $I_5$ holds similarly as the second item above after the first "$\lesssim$". Meanwhile, estimates of $I_6$ and $I_7$ hold similarly with that of $I_1$ and $I_2$ respectively. Combining those estimates above in this section, \eqref{2.2} holds for any $r\in[R,\,2R]$, i.e.
\be
\left|\int_{-\pi}^\pi\int_0^{2\pi} p(r,\th,z)-p(R,\th,z)d\th dz\right|\lesssim 1,\quad\forall r\in[R,\,2R].
\ee
Applying the mean value theorem, for a fixed $R>1$ and $r\in[R,\,2R]$, there exist $\theta(r)\in[0,\,2\pi]$ and $z(r)\in[-\pi,\,\pi]$, such that
\be
\big|p(r,\th(r),z(r))-p(R,\th(r),z(r))\big|=\left|\int_{-\pi}^\pi\int_0^{2\pi} p(r,\th,z)-p(R,\th,z)d\th dz\right|\lesssim 1.
\ee
Combining this with the uniformly boundedness of $\p_zp$ and $\p_\th p$, it follows that, $\forall\th\in[0,\,2\pi]$, $z\in[-\pi,\,\pi]$
\be
\begin{split}
|p(r,\th,z)-p(R,\th,z)|\leq&|p(r,\th,z)-p(r,\th,z(r))|+|p(r,\th,z(r))-p(r,\th(r),z(r))|\\
&+|p(r,\th(r),z(r))-p(R,\th(r),z(r)|\\
&+|p(R,\th(r),z(r))-p(R,\th(r),0)|\\
&+|p(R,\th(r),0)-p(R,0,0)|\\
\leq&\,\,|\p_zp|\cdot|z-z(r)|+|\p_\th p|\cdot|\th-\th(r)|+C\\
&+|\p_zp|\cdot|z(r)|+|\p_\th p|\cdot|\th(r)|\\
\lesssim& 1.
\end{split}
\ee
Hence we have
\be\label{2.12}
\sup_{r\in[R,\,2R],\,\th\in[0,\,2\pi],\,z\in[-\pi,\,\pi]}\left|p(r,\th,z)-p(R,0,0)\right|\lesssim1.
\ee

\qed

\subsection{Trivialness of $u$ and $h$}\label{SEC2.4}
At the beginning, we claim that $u^r,\,h^r\in L^{2}(\mathbb{R}^2\times\mathbb{T})$. The reason is: according to the divergence-free condition and $u^\th\equiv0$, we see that
\be\label{3.1}
\p_r(ru^r(r,\th,z))+\p_z(ru^z(r,\th,z))=0.
\ee
Integrating  \eqref{3.1} on $z$ from $-\pi$ to $\pi$, it follows that
\be
\p_r\left(r\int_{-\pi}^\pi u^r(r,\th,z)dz\right)=-\int_{-\pi}^\pi\p_z(ru^z(r,\th,z))dz=-ru^z(r,\th,z)\Big|_{z=-\pi}^\pi=0.
\ee
Here the last identity follows from the periodic condition of $u$ in $z$-direction. This leads to
\be
r\int_{-\pi}^\pi u^r(r,\th,z)dz=C(\th),
\ee
where $C(\th)$ is a function depends only on $\th$. Moreover, we find $C(\th)\equiv0$ by choosing $r=0$. Therefore
\be
\int_{-\pi}^\pi u^r(r,\th,z)dz=0.
\ee
Hence we have, by using the Poincar\'{e} inequality and the D-solution condition
\be
\begin{split}
\int_{-\pi}^\pi\int_{\mathbb{R}^2} |u^r|^2dx=&\int_{\mathbb{R}^2}\int_{-\pi}^\pi\left|u^{r}(x',z)-\frac{1}{2\pi}\int_{-\pi}^\pi u^r(x',z')dz'\right|^2dzdx'\\
\lesssim&\int_{\mathbb{R}^2}\int_{-\pi}^\pi|\p_zu^r(x',z)|^2dzdx'\\
\leq&\int_{\mathbb{R}^2\times\mathbb{T}}|\nabla u(x)|^2dx<\infty.
\end{split}
\ee
At the same time, the related estimate holds for $h^r$, that is
\be
\int_{-\pi}^\pi\int_{\mathbb{R}^2} |h^r|^2dx<\infty,
\ee
since $h$ is divergence-free and axially symmetric, which means equation \eqref{3.1} also holds for $h$. The rest is similar with that of $u$ and the claim is proved.

Now let $\phi=\phi(\rho)$ be a smooth cut-off function satisfying
\be
\left\{
\begin{aligned}
\phi(\rho)&=1,\quad\rho\in[0,1],\\
\phi(\rho)&=0,\quad\rho\geq2,\\
0\leq&\phi\leq1,\quad\forall\rho\in[0,\infty),\\
\end{aligned}
\right.
\ee
with $\phi'$ and $\phi''$ being bounded. And we set $\phi_R(y')=\phi\left(\frac{|y'|}{R}\right)$ with $y'\in\mathbb{R}^2$ and $R>0$. Testing the first equation of \eqref{1.3}
\be
u\cdot\nabla u+\nabla p-h\cdot\nabla h-\Delta u=0
\ee
with $u\phi_R$, we achieve that
\be
\int_{\mathbb{R}^2\times\mathbb{T}}u\phi_R\Delta udx=\int_{\mathbb{R}^2\times\mathbb{T}}u\phi_R\left(u\cdot\nabla u-h\cdot\nabla h+\nabla(p-p(R,0,0))\right)dx.
\ee
Direct integrating by parts implies
\be\label{2.25}
\begin{split}
\int_{\mathbb{R}^2\times\mathbb{T}}|\nabla u|^2\phi_Rdx-&\frac{1}{2}\int_{\mathbb{R}^2\times\mathbb{T}}|u|^2\Delta\phi_Rdx\\
=&\frac{1}{2}\int_{\mathbb{R}^2\times\mathbb{T}}|u|^2u\cdot\nabla\phi_Rdx+\int_{\mathbb{R}^2\times\mathbb{T}}(p(r,\th,z)-p(R,0,0))u\cdot\nabla\phi_Rdx\\
&-\sum_{i,j=1}^3\int_{\mathbb{R}^2\times\mathbb{T}}h_ih_j\p_{x_i}u_j\phi_Rdx-\sum_{i,j=1}^3\int_{\mathbb{R}^2\times\mathbb{T}}h_ih_ju_j\p_{x_i}\phi_Rdx
\end{split}
\ee
Meanwhile, we test the second equation of \eqref{1.3}
\be
u\cdot\nabla h-h\cdot\nabla u-\Delta h=0
\ee
with $h\phi_R$ to get
\be\label{2.27}
\int_{\mathbb{R}^2\times\mathbb{T}}h\phi_R\Delta hdx=\int_{\mathbb{R}^2\times\mathbb{T}}h\phi_R\left(u\cdot\nabla h-h\cdot\nabla u\right)dx.
\ee
Integrating by parts, \eqref{2.27} is equivalent to
\be\label{2.28}
\begin{split}
\int_{\mathbb{R}^2\times\mathbb{T}}|\nabla h|^2\phi_Rdx-&\frac{1}{2}\int_{\mathbb{R}^2\times\mathbb{T}}|h|^2\Delta\phi_Rdx\\
=&\frac{1}{2}\int_{\mathbb{R}^2\times\mathbb{T}}|h|^2u\cdot\nabla\phi_Rdx+\sum_{i,j=1}^3\int_{\mathbb{R}^2\times\mathbb{T}}h_ih_j\p_{x_i}u_j\phi_Rdx
\end{split}
\ee
Therefore, the following equation is achieved by adding \eqref{2.25} and \eqref{2.28} together:
\be
\begin{split}
\int_{\mathbb{R}^2\times\mathbb{T}}\left(|\nabla u|^2+|\nabla h|^2\right)&\phi_Rdx-\frac{1}{2}\int_{\mathbb{R}^2\times\mathbb{T}}\left(|u|^2+|h|^2\right)\Delta\phi_Rdx\\
=&\int_{\mathbb{R}^2\times\mathbb{T}}\left(\frac{1}{2}|u|^2+\frac{1}{2}|h|^2+(p(r,\th,z)-p(R,0,0))\right)u\cdot\nabla\phi_Rdx\\
&-\int_{\mathbb{R}^2\times\mathbb{T}}(h\cdot u)(h\cdot\nabla\phi_R)dx.
\end{split}
\ee
We denote $\bar{B}_{2R/R}:=\left\{x':\,R\leq|x'|\leq2R\right\}$ the dyadic annulus. Since $\phi_R$ depends only on $r$, it follows that
\be\label{2.31}
\begin{split}
&\int_{\mathbb{R}^2\times\mathbb{T}}\left(|\nabla u|^2+|\nabla h|^2\right)\phi_Rdx\\
\leq&\int_{-\pi}^\pi\int_{\bar{B}_{2R/R}}\left(|u|^2+|h|^2\right)\cdot|\Delta\phi_R|dx'dz\\
&+\int_{-\pi}^\pi\int_{\bar{B}_{2R/R}}|u^r|\cdot|\p_r\phi_R|\cdot(|u|^2+|h|^2)dx'dz\\
&+\sup_{r\in[R,\,2R],\,\th\in[0,\,2\pi],\,z\in[-\pi,\,\pi]}\left|p(r,\th,z)-p(R,0,0)\right|\int_{-\pi}^\pi\int_{\bar{B}_{2R/R}}|u^r|\cdot|\p_r\phi_R|dx'dz\\
&+\int_{-\pi}^\pi\int_{\bar{B}_{2R/R}}|h|\cdot|u|\cdot|h^r|\cdot|\p_r\phi_R|dx'dz\\
:=&I_1+I_2+I_3+I_4.
\end{split}
\ee
First, $I_1$ satisfies
\be
\begin{split}
I_1&\lesssim\frac{1}{R^2}\int_{-\pi}^\pi\int_{\bar{B}_{2R/R}}(|u|^2+|h|^2)dx'dz\\
&\leq\frac{1}{R^2}\left(\|u\|^2_{L^\infty(\bar{B}_{2R/R})}+\|h\|^2_{L^\infty(\bar{B}_{2R/R})}\right)\int_{-\pi}^\pi\int_{\bar{B}_{2R/R}}dx'dz\\
&\lesssim\|u\|^2_{L^\infty(\bar{B}_{2R/R})}+\|h\|^2_{L^\infty(\bar{B}_{2R/R})}\to 0,\quad\text{ as }\quad R\to\infty.
\end{split}
\ee
Here we applied the vanishing of both $u$ and $h$ at the far field, and the same as below for the estimates of $I_2$ and $I_4$. Using the H\"{o}lder inequality, $I_2$ follows that
\be
\begin{split}
I_2&\lesssim\frac{1}{R}\int_{-\pi}^\pi\int_{\bar{B}_{2R/R}}|u^r|\cdot\left(|u|^2+|h|^2\right)dx'dz\\
&\leq\frac{1}{R}\left(\|u\|^2_{L^\infty(\bar{B}_{2R/R})}+\|h\|^2_{L^\infty(\bar{B}_{2R/R})}\right)\left(\int_{-\pi}^\pi\int_{\bar{B}_{2R/R}}|u^r|^2dx'dz\right)^{1/2}\cdot|\bar{B}_{2R/R}|^{1/2}\\
&\lesssim\left(\|u\|^2_{L^\infty(\bar{B}_{2R/R})}+\|h\|^2_{L^\infty(\bar{B}_{2R/R})}\right)\cdot\|u^r\|_{L^2(\mathbb{R}^2\times\mathbb{T})}\to0,\quad\text{ as }\quad R\to\infty.
\end{split}
\ee
Next, for $I_3$ we have
\be
\begin{split}
I_3\lesssim&\sup_{r\in[R,\,2R],\,\th\in[0,\,2\pi],\,z\in[-\pi,\,\pi]}\left|p(r,\th,z)-p(R,0,0)\right|\\
&\hskip 3cm\cdot\frac{1}{R}\left(\int_{-\pi}^\pi\int_{\bar{B}_{2R/R}}|u^r|^2dx'dz\right)^{1/2}\cdot|\bar{B}_{2R/R}|^{1/2}\\
\lesssim&\|u^r\|_{L^2(\bar{B}_{2R/R})}\to0,\quad\text{ as }\quad R\to\infty.
\end{split}
\ee
Here we have applied the H\"older inequality and the boundedness of the oscillation of $p$ in dyadic annulus which is achieved in Lemma \ref{Lemma2.2}. Finally, the following estimate is satisfied by $I_4$
\be
\begin{split}
I_4&\lesssim\frac{1}{R}\cdot\|u\|_{L^\infty(\bar{B}_{2R/R})}\cdot\|h\|_{L^\infty(\bar{B}_{2R/R})}\cdot\left(\int_{-\pi}^\pi\int_{\bar{B}_{2R/R}}|h^r|^2dx'dz\right)^{1/2}\cdot|\bar{B}_{2R/R}|^{1/2}\\
&\lesssim\|u\|_{L^\infty(\bar{B}_{2R/R})}\cdot\|h\|_{L^\infty(\bar{B}_{2R/R})}\cdot\|h^r\|_{L^2(\bar{B}_{2R/R})}\to0,\quad\text{ as }\quad R\to\infty.
\end{split}
\ee
Combining those estimates of $I_1$, $I_2$, $I_3$ and $I_4$, \eqref{2.31} implies
\be
\int_{\mathbb{R}^2\times\mathbb{T}}\left(|\nabla u|^2+|\nabla h|^2\right)dx=0,
\ee
by choosing $R\to\infty$. This means $u$ and $h$ are both constants. Recalling $u$ and $h$ vanish at the far field, we deduce the trivialness of $u$ and $h$ themselves. Now we have finished the proof of the part (i) of \textbf{Theorem \ref{thm1.1}}.

\qed

\section{Proof of Corollary 1.1, part (i)}
This section is devoted to the case that $u^\th\equiv0$ and $h$ is axially symmetric. First we see, in this situation, \eqref{1.2} turns to
\begin{equation}\label{1.21}
\left\{
\begin{aligned}
&(u^r\p_r+u^z\p_z)u^r+\p_r p=(h^r\p_r+h^z\p_z)h^r-\frac{(h^\th)^2}{r}+(\Delta-\frac{1}{r^2})u^r, \\
&-\frac{2}{r^2}\p_\th u^r+\frac{1}{r}\p_\th p=(h^r\p_r+h^z\p_z)h^\th+\frac{h^rh^\th}{r}, \\
&(u^r\p_r+u^z\p_z)u^z+\p_z p=(h^r\p_r+h^z\p_z)h^z+\Delta u^z ,                                    \\
&(u^r\p_r+u^z\p_z)h^r-(h^r\p_r+\frac{1}{r}h^\th\p_\th+h^z\p_z)u^r=(\Delta-\frac{1}{r^2})h^r,\\
&(u^r\p_r+u^z\p_z)h^\th-\frac{u^rh^\th}{r}=(\Delta-\frac{1}{r^2})h^\th,\\
&(u^r\p_r+u^z\p_z)h^z-(h^r\p_r+\frac{1}{r}h^\th\p_\th+h^z\p_z)u^z=\Delta h^z,\\
&\nabla\cdot u=\p_ru^r+\frac{u^r}{r}+\p_zu^z=0,\quad \nabla\cdot h=\p_rh^r+\frac{h^r}{r}+\p_zh^z=0.
\end{aligned}
\right.
\end{equation}
We point out that, in this situation, $h^\th$ is vanishing so that the proof follows from the proof of Theorem \ref{thm1.1}. Here goes the proof of the vanishing of $h^\th$.
\subsection*{The Vanishing of $h^\th$}
Under the axially symmetric condition of $h$ and the vanishing of $u^\th$, the equation of $h^\th$ reads
\be
(u^r\p_r+u^z\p_z)h^\th-\frac{u^rh^\th}{r}=(\Delta-\frac{1}{r^2})h^\th.
\ee
Denoting $H=\frac{h^\th}{r}$, direct calculation shows that $H$ satisfies
\be\label{EQH}
\left(\Delta+\frac{2}{r}\p_r\right)H-(u^r\p_r+u^z\p_z)H=0.
\ee
Since $h$ is axially symmetric, the Laplacian operator here write
\be
\Delta=\frac{\p^2}{\p r^2}+\frac{1}{r}\frac{\p}{\p r}+\frac{\p^2}{\p z^2}.
\ee
Therefore, if we denoting
\be
\Delta_5:=\sum_{i=1}^4\frac{\p^2}{\p x_i^2}+\frac{\p^2}{\p z^2}
\ee
and $r=\sqrt{x_1^2+x_2^2+x_3^2+x_4^2}$, \eqref{EQH} becomes
\be
\Delta_5H-(u^r\p_r+u^z\p_z)H=0.
\ee
See \cite{Koch2009} or \cite{LiPan2019} for more details about this "dimension lifting method". From the boundedness of $h^\th$, one find
\be
\lim_{r\to\infty}H=0
\ee
uniformly for all $z$. Therefore, $H\equiv0$ is achieved by the maximum princlple. This leads to the vanishing of $h^\th$.
\section{Proof of Corollary 1.1, part (ii)}
In this section, we consider the case that both $u$ and $h$ are axially symmetric. At the beginning we see \eqref{1.2} now turns to
\begin{equation}\label{1.22}
\left\{
\begin{aligned}
&(u^r\p_r+u^z\p_z)u^r -\frac{(u^\th)^2}{r}+\p_r p=(h^r\p_r+h^z\p_z)h^r-\frac{(h^\th)^2}{r}+(\Delta-\frac{1}{r^2})u^r, \\
&(u^r\p_r+u^z\p_z) u^\th+\frac{u^r u^\th}{r}=(h^r\p_r+h^z\p_z)h^\th+\frac{h^rh^\th}{r}+(\Delta-\frac{1}{r^2})u^\th , \\
&(u^r\p_r+u^z\p_z)u^z+\p_z p=(h^r\p_r+h^z\p_z)h^z+\Delta u^z ,                                    \\
&(u^r\p_r+u^z\p_z)h^r-(h^r\p_r+h^z\p_z)u^r=(\Delta-\frac{1}{r^2})h^r,\\
&(u^r\p_r+u^z\p_z)h^\th-(h^r\p_r+h^z\p_z)u^\th+\frac{u^\th h^r}{r}-\frac{h^\th u^r}{r}=(\Delta-\frac{1}{r^2})h^\th,\\
&(u^r\p_r+u^z\p_z)h^z-(h^r\p_r+h^z\p_z)u^z=\Delta h^z,\\
&\nabla\cdot u=\p_ru^r+\frac{u^r}{r}+\p_zu^z=0,\quad \nabla\cdot h=\p_rh^r+\frac{h^r}{r}+\p_zh^z=0.
\end{aligned}
\right.
\end{equation}
The main idea are similar with the proof of Theorem \ref{thm1.1} so that we only focus on the different portion. First we see, by Section \ref{SEC2.1}, $u$ and $h$ are bounded up to their second order derivatives. Combining these with the third equation of \eqref{1.22}, we have the boundedness of $\p_zp$. Integrating the first equation in \eqref{1.22} like Section 2, we are ready to prove the boundedness of $p$ in dyadic annulus. Pay attention that, due to the axially symmetric condition for both $u$ and $h$, $p$ is no longer a function of $\th$. We only need to prove the boundedness of
\be
I:=\left|\int_R^{r_0}\int_{-\pi}^\pi\int_0^{2\pi}\frac{(u^\th)^2}{r}d\th dzdr\right|+\left|\int_R^{r_0}\int_{-\pi}^\pi\int_0^{2\pi}\frac{(h^\th)^2}{r}d\th dzdr\right|
\ee
since the boundedness of the rest terms have already proven in Section 2. Here goes the boundedness of $I$:
\be
\begin{split}
I\lesssim&\left(\|u^\th\|_{L^\infty(\bar{B}_{2R/R})}^2+\|h^\th\|_{L^\infty(\bar{B}_{2R/R})}^2\right)\int_{R}^{2R}\frac{1}{r}dr\\
\lesssim&\|u^\th\|_{L^\infty(\bar{B}_{2R/R})}^2+\|h^\th\|_{L^\infty(\bar{B}_{2R/R})}^2\\
\lesssim&1.
\end{split}
\ee
Then the vanishing of $u$ and $h$ is achieved by following the method in Section \ref{SEC2.4}. We omit the details here.

\section{Appendix: Some details of the boundary conditions}

This Appendix is devoted to some explanations of the boundary conditions in Corollary \ref{Cor2}. As we mentioned in the Corollary \ref{Cor2}, instead of the periodic condition for the velocity field, our method is also valid for a certain Navier slip boundary condition with a slight modification. That is
\be\label{NavierB}
u\cdot n=0,\quad\left(\mathbb{D}u\cdot n\right)_{\tau}=0,\quad\forall x\in\p\Omega.
\ee
Here $n$ is the outward unit normal to $\Omega$. $\mathbb{D}$ is the strain tensor
\be
\mathbb{D}u=\frac{1}{2}\left(\nabla u+\nabla^{T}u\right).
\ee
And for a vector filed $v$, $v_\tau$ stands for its tangential part: $v_\tau=v-(v\cdot n)n$. In our case, since $\Omega=\mathbb{R}^2\times[-\pi,\,\pi]$, we have $n=(0,0,\pm1)$. Therefore, \eqref{NavierB} is reduced to
\be\label{NavierB1}
u^z=0,\quad\p_zu_1=0,\quad\p_zu_2=0,\quad\forall z=-\pi\text{ or }\pi.
\ee
In the cylinder coordinate, \eqref{NavierB1} equals to
\be
\left\{
\begin{aligned}
&u^z=0,\\
&\p_zu^r\cos\th-\p_zu^\th\sin\th=0,\\
&\p_ru^z\sin\th+\p_zu^\th\cos\th=0,
\end{aligned}
\right.
\quad\quad\quad\forall z=-\pi\text{ or }\pi.
\ee
That is,
\be\label{2.45}
\p_zu^r\Big|_{z=-\pi,\,\pi}=\p_zu^\th\Big|_{z=-\pi,\,\pi}=u^z\Big|_{z=-\pi,\,\pi}\equiv0.
\ee
Meanwhile, for magnetic field $h$, our method is valid for the Dirichlet condition
\be
h=0,\quad\forall z=-\pi\text{ or }\pi,
\ee
and the following two physical conditions, which are widely used in the research of the boundary value problem or the initial-boundary value problems to the MHD system. See \cite{Huang-Wang2012, Wu2011}, etc.
\be
\text{[PC1]}
\left\{
\begin{aligned}
&h\cdot n=0,\\
&\nabla\times h\times n=0,\\
\end{aligned}
\right.
\quad\quad
\text{[PC2]}
\left\{
\begin{aligned}
&h\cdot n=0,\\
&\nabla\times h=0,\\
\end{aligned}
\right.
\quad\quad\forall x\in\p\Omega.
\ee
In our cases, similarly to \eqref{2.45} before, \emph{[PC1]} and \emph{[PC2]} are equivalent and both of them can be simplified to
\be
\p_zh^r\Big|_{z=-\pi,\,\pi}=\p_zh^\th\Big|_{z=-\pi,\,\pi}=h^z\Big|_{z=-\pi,\,\pi}\equiv0.
\ee

\section*{Acknowledgments} The authors wish to thank Prof. Huicheng Yin in Nanjing Normal University and Prof. Qi S. Zhang in UC Riverside for their constant encouragement on this topic.

\end{document}